.
.
.
\font\script=eusm10.
\font\sets=msbm10.
\font\stampatello=cmcsc10.
.
\def\0{{\bf 0}}
\def\1{{\bf 1}}
\def\defineq{\buildrel{def}\over{=}}

\def\definiz{\buildrel{def}\over{\Longleftrightarrow}}
\def\C{\hbox{\sets C}}

\def\N{\hbox{\sets N}}

\def\Primes{\hbox{\sets P}}

\def\Z{\hbox{\sets Z}}

\def\square{\hbox{\vrule\vbox{\hrule\phantom{s}\hrule}\vrule}}

\def\IrrP1F{{\hbox{\rm Irr}^{(P)}_1\,F}}

\def\Hil{{\rm Hil}}

\par
\centerline{\bf Absolute convergence of Ramanujan expansions}
\par
\centerline{\bf admits coefficients' coexistence of Ramanujan expansions}
\bigskip
\centerline{Giovanni Coppola}\footnote{ }{MSC $2020$: $11{\rm N}37$ - Keywords: arithmetic function, Ramanujan expansion, absolute convergence} 

\bigskip

\par
\noindent
{\bf Abstract}. In this self-contained short note, we prove that {\it every arithmetic function $F:\N\rightarrow\C$ has infinitely many Ramanujan coefficients $G:\N\rightarrow\C$ giving an absolutely convergent Ramanujan expansion for $F$}. This is \lq \lq coefficients' coexistence\rq \rq: the non-uniqueness, once fixed any $F$, of these $G$. They are infinitely many ! 

\bigskip

\par
\noindent{\bf 1. Introduction. Notations and well-known properties} 
\bigskip
\par
\noindent
The notation we need soon is for sets of numbers: $\N\defineq \{1,2,\ldots\}$ are the {\it naturals}, i.e., positive integers; while we denote $\N_0\defineq \{0,1,2,\ldots\}$ {\it non-negative integers}; as usual, $\C$ is the set of {\it complex numbers}. In the following, for integers $d,n$ with $d\neq 0$ we abbreviate $d|n$ for: $d$ {\it divides} $n$, in other words $n$ is a multiple of $d$. 
\par
\noindent
The set of {\it prime numbers} is abbreviated $\Primes$, while they are denoted as $p$, with or without subscripts. Two integers are called {\it coprime}, by definition, IFF (If and only If) there's no $p$ dividing both: in other words, their {\it greatest common divisor}, g.c.d., is $1$. When not explicitly stated otherwise, $n\le x$ means $n\ge 1$, too. 

\bigskip

\par
\noindent
Ramanujan studied, in the history-making paper [R], the {\stampatello Ramanujan Sum} with {\it modulus} $q$ \& {\it argument} $a$ :
$$
\forall q\in \N,\enspace \forall a\in \Z,
\qquad
c_q(a)\defineq \sum_{{j\le q}\atop {(j,q)=1}}\cos {{2\pi ja}\over q}, 
$$
\par
\noindent
where $(j,q)$ abbreviates the g.c.d. of $j$ and $q$. See that $c_q(0)=\varphi(q)\defineq |\{n\le q : (n,q)=1\}|$ is {\it Euler's Totient}; while $c_q(1)=\mu(q)$ is {\it M\"obius function}, that is $\mu(q)\defineq 0$ on $q\in \N$ multiples of the square of some prime $p$, say $p^2|q$, instead on {\it square-free} $q\in \N$ (no prime $p$ has $p^2|q$) it's $\mu(q)\defineq (-1)^{\omega(q)}$, where $\omega(q)\defineq |\{p\in \Primes : p|q\}|$ is the {\it number of prime-factors} of our $q\in \N$. In particular, $\omega(1)=0$ entails $\mu(1)=1$.
\par
These two instances are particular cases of H\"older's Formula [M], compare Chapter 26 [D] : 
$$
c_q(a)=\varphi(q){{\mu(q/(a,q))}\over {\varphi(q/(a,q))}},
\quad
\forall q\in \N, \forall a\in \Z
\leqno{(1)}
$$
\par
\noindent
that, in particular, ensures that {\it all Ramanujan Sums $c_q(a)$ are integers}, whatever $a,q$ above. 
\par
Surprisingly enough, we may use them as a kind of \lq \lq {\it harmonics}\rq \rq, say, for the {\it Arithmetic Functions}; these are all the functions 
$$
F:\N \rightarrow \C
$$
\par
\noindent
and this is, actually, the idea of Ramanujan [R] : given $F$, we look for another arithmetic function $G:\N \rightarrow \C$ that we call a {\it Ramanujan coefficient} for $F$, by definition, IFF following limit exists, over $x\to \infty$, in $\C$ and 
$$
\forall a\in \N,
\qquad
F(a)=\sum_{q=1}^{\infty}G(q)c_q(a)
 \defineq \lim_x \sum_{q\le x}G(q)c_q(a)\in \C, 
$$
\par
\noindent
a formula that, in fact, is called the {\it Ramanujan Expansion for $F$} with coefficient $G$. Ramanujan Himself [R] gave, for the {\it null-function}, defined as 
$$
\0(a)\defineq 0,
\qquad
\forall a\in \N,
$$
\par
\noindent
the Ramanujan Expansion 
$$
\forall a\in \N,
\qquad
\0(a)=\sum_{q=1}^{\infty}{1\over q}c_q(a). 
$$
\par				
\noindent
The Proof of this formula, in spite of its apparent simplicity, involves a deep property of prime numbers, namely the PNT, {\it Prime Number Theorem} [M],[R] (compare [C2], too). 
\par
For present aims, it's more important to focus on two other important properties coming from previous null-function expansion. First, it proves that, at the beginning for $F=\0$, there is {\it no uniqueness of the Ramanujan coefficient} $G$ (the other Ramanujan coefficient is the trivial $G=\0$, for $\0$); second, the series {\it doesn't converge absolutely}: at least for $a=1$, in which previous series is $\sum_{q=1}^{\infty}\mu(q)/q$ and converges to $0$ (compare \lq \lq M\"obius Language PNT\rq \rq, in [C2] Fact 1.6), BUT, $\sum_{q=1}^{\infty}\mu^2(q)/q=+\infty$, because $\mu^2(q)=1$ IFF $q$ is square-free and we have \enspace $\sum_{p=2}^{\infty}1/p=+\infty$, i.e., the divergence of primes' reciprocals series, a very well-known property [T]. Also, non-uniqueness of Ramanujan coefficients $G$ holds {\it for all} $F:\N \rightarrow \C$. Actually, this follows as a particular consequence of our next Corollary.  
\smallskip
\par
\noindent
These two properties, considered together, may seem to suggest that all Ramanujan coefficients $G\neq \0$ for the null-function give a non-absolutely convergent Ramanujan expansion... But this is dramatically false ! 
\smallskip
\par
\noindent
In our next Theorem, the null-function $\0$ reveals, so to speak: WE HAVE AN INFINITY OF RAMANUJAN COEFFICIENTS WITH ABSOLUTELY CONVERGENT RAMANUJAN EXPANSIONS for $\0$. Also, as a Corollary, we get the same, more in general, FOR ANY FIXED ARITHMETIC FUNCTION. See next $\S2$. 
\medskip
\par
Notice that, for general $F$, writing
$$
\forall a\in \N, 
\qquad
\sum_{q=1}^{\infty}G(q)c_q(a)=F(a)
$$
\par
\noindent
means : $\forall a\in \N$ the series converges and its sum is $F(a)$. Furthermore, the {\stampatello absolute convergence} of this series, $\forall a\in\N$ fixed, {\stampatello ensures both its convergence \& its sum-independence}, $\forall a\in\N$ fixed, {\stampatello whatever summation method} is {\stampatello chosen.} (Not only the classic one, $q\le x$, then $x\to \infty$, we saw above.)
\par
In passing, see that \lq \lq {\stampatello $G$ is a Ramanujan Coefficient for $F$\rq \rq \enspace entails} that the corresponding {\stampatello Ramanujan Expansion for $F(a)$ holds $\forall a\in \N$, as above}. (Esp., Definition 1.1 in [L] hasn't such a limitation.) 
\medskip
\par
\noindent
Another, say, VIP, Very Important Property, we need uses the classic definition of {\it $p-$adic valuation} of a natural number $q$: 
$$
v_p(q)\defineq \max\{K\in \N_0 : p^K|q\}
$$
\par
\noindent
and it comes from $(1)$ very easily: 
$$
c_q(a)\neq 0
\quad \Rightarrow \quad
v_p(q)\le v_p(a)+1,\enspace \forall p|q
\leqno{(2)}
$$
\par
\noindent
a kind of prime-powers limit for non-zero Ramanujan Sums, which we call {\stampatello Ramanujan Vertical Limit}. Actually, also converse implication (not needed here) holds in $(2)$, as is proved in [C3], Proposition 3.1. 
\par
When coming to bounds on absolute values, in Ramanujan expansions, of course following bound is very useful: 
$$
|c_q(a)|\leq (q,a)
\leq a, 
\quad 
\forall q\in \N,\enspace \forall a\in \N,
\leqno{(3)}
$$
\par
\noindent
where the first inequality is Lemma 2.1 [CM] (see there, for the elementary proof, again using $(1)$ above), but here we only need second inequality (coming from first one, trivially because the g.c.d of two naturals doesn't go beyond both of them). 

\bigskip

\centerline{\bf Paper's Plan}
\smallskip
\item{$\star$} in next section, $\S2$, we give our results with their Proofs, starting from above numbered formul\ae; 
\item{$\star$} in $\S3$ we give complementary, deeper and noteworthy aspects, so to speak, of our results; 
\item{$\star$} in the Appendix, we offer a particular Proof \& details not strictly needed for our results. 
 
\vfill
\eject

\par				
\noindent{\bf 2. Results: Statements and Proofs} 
\bigskip
\par
\noindent
We recall the definition of {\it multiplicative} arithmetic function $F:\N \rightarrow \C$, namely:
\smallskip
\par
\centerline{$a$ and $b$ coprime entails $F(ab)=F(a)F(b)$.}
\smallskip
\par
We don't, for the time being (compare our other papers, too), assume that $F(1)=1$; in this way, we may have $F(1)=0$, so to include $\0$ in the multiplicative functions. Other Authors, instead, assume $F(1)\neq 0$ which, actually, is equivalent to $F(1)=1$. 
\par
\noindent
Anyway, a multiplicative arithmetic function is clearly individuated by its values on prime-powers: compare next Lemma Proof. By the way, we recall: \enspace once fixed $a\in \Z$, $c_q(a)$ is a multiplicative function of $q\in \N$, see quoted [D]. 
\bigskip
\par
We start with a property that we found in [CG1], compare the \lq \lq {\it exotic}\rq \rq \enspace {\it Ramanujan coefficients} for the null-function. We prove it here in a completely self-contained way. \enspace See $\S3$ for further details. 
\smallskip
\par
\noindent{\bf Lemma}. {\stampatello (Null-function's Ramanujan coefficients with absolutely convergent expansions)}
\par
\noindent
{\it Let } $G:\N \rightarrow \C$ {\it be multiplicative and assume that } $\exists p_0\in \Primes$, {\it with}
$$
G(p_0^K)=1,
\qquad
\forall K\in \N_0.
$$
\par
\noindent
{\it Furthermore, let the series with indices coprime to $p_0$ , say}
$$
\sum_{{q=1}\atop {(q,p_0)=1}}^{\infty}G(q)c_q(a),
\qquad
\hbox{\it converge\enspace absolutely}
\enspace
\forall a\in \N\enspace \hbox{\it fixed}.
$$
\par
\noindent
{\it Then}
$$
\forall a\in \N,
\qquad
\sum_{q=1}^{\infty}\left|G(q)c_q(a)\right|<\infty
$$
\par
\noindent
{\it and}
$$
\forall a\in \N,
\qquad
\sum_{q=1}^{\infty}G(q)c_q(a)=0;
$$
\par
\noindent
{\it these two conditions simply say that $\0(a)$ has an absolutely convergent Ramanujan expansion, converging $\forall a\in \N$ fixed, of } {\stampatello Ramanujan coefficient} $G:\N \rightarrow \C$. 
\smallskip
\par
\noindent{\bf Proof}. First of all, fix $a\in \N$; writing any $q\in \N$ as $q=p_0^K\,m$, $K\in \N_0$ and $(m,p_0)=1$ : 
$$
\sum_{q=1}^{\infty}\left|G(q)c_q(a)\right|
=\sum_{0\le K\le v_{p_0}(a)+1}\left|c_{p_0^K}(a)\right|\sum_{{m=1}\atop {(m,p_0)=1}}^{\infty}\left|G(m)c_m(a)\right|
 =\left(\sum_{K=0}^{v_{p_0}(a)+1}\left|c_{p_0^K}(a)\right|\right)\sum_{{q=1}\atop {(q,p_0)=1}}^{\infty}\left|G(q)c_q(a)\right|
  <\infty, 
$$
\par
\noindent
because, fixed $a\in \N$, $|G(q)c_q(a)|$ is {\stampatello multiplicative} w.r.t. $q\in \N$, $c_{p_0^K}(a)=0$ $\forall K>v_{p_0}(a)+1$ by {\stampatello Ramanujan vertical limit} $(2)$ and we use our absolute convergence hypothesis on $(q,p_0)=1$ indices. As we know, in particular, that the LHS (Left Hand Side) series converges $\forall a\in \N$ {\stampatello fixed}, we get : 
$$
\sum_{q=1}^{\infty}G(q)c_q(a)
 =\left(\sum_{K=0}^{v_{p_0}(a)+1}c_{p_0^K}(a)\right)\sum_{{q=1}\atop {(q,p_0)=1}}^{\infty}G(q)c_q(a)=0, 
$$
\par
\noindent
again by previous properties \& hypothesis, where now we also use from H\"older's Formula $(1)$: $\forall a\in \N$ {\stampatello fixed},
$$
\sum_{K=0}^{v_{p_0}(a)+1}c_{p_0^K}(a)=\sum_{K=0}^{v_{p_0}(a)}\varphi(p_0^K)-p_0^{v_{p_0}(a)}
 =1+\sum_{K\le v_{p_0}(a)}\left(p_0^K-p_0^{K-1}\right)-p_0^{v_{p_0}(a)}
  =0.
$$
\rightline{\square}

\vfill
\eject

\par				
For the {\it Ramanujan clouds} of a fixed $F$, i.e., particular subsets of Ramanujan coefficients for $F$, see our papers [C1], [CG1] and [CG2]; here we only need the following.
\smallskip
\par
\noindent{\bf Definition}. {\stampatello (Absolute Ramanujan Clouds)}
\par
\noindent
{\it For the} {\stampatello null-function} $\0:\N\rightarrow \C$ {\it we introduce the set}
$$
<\0>_{\rm abs}\defineq \left\{ G:\N\rightarrow \C \;\left|\; \forall a\in \N, \sum_{q=1}^{\infty}G(q)c_q(a)=0\enspace \hbox{\stampatello and}\enspace \sum_{q=1}^{\infty}|G(q)c_q(a)|<\infty \right. \right\}, 
$$
\par
\noindent
{\it called the} {\stampatello Absolute Ramanujan Cloud} {\it of} $\0$, {\stampatello i.e.,} {\it the set of all the} {\stampatello Ramanujan Coefficients} {\it of} $\0$ {\stampatello with} {\it an} {\stampatello absolutely convergent Ramanujan Expansion} {\it of} $\0$.\thinspace(Trivially, $\0\in <\0>_{\rm abs}$, so it's not empty.) 
\par
{\it More in general,}
\par
\noindent
{\stampatello All} ({\it fixed}) {\stampatello arithmetic functions} $F$ {\it have their} {\stampatello Absolute Ramanujan Cloud}: 
$$
<F>_{\rm abs}\defineq \left\{ G:\N\rightarrow \C \;\left|\; \forall a\in \N, \sum_{q=1}^{\infty}G(q)c_q(a)=F(a)\enspace \hbox{\stampatello and}\enspace \sum_{q=1}^{\infty}|G(q)c_q(a)|<\infty \right. \right\} 
$$
\par
\noindent
{\it as the set of all the} {\stampatello Ramanujan Coefficients} {\it for} $F$ {\stampatello with} {\it an} {\stampatello absolutely convergent Ramanujan Expansion} {\it for} $F$. (Next Proposition entails that any fixed $F:\N \rightarrow \C$ has non-empty $<F>_{\rm abs}$.)

\medskip

Our Definition allows us to abbreviate next result's statement.
\smallskip
\par
\noindent{\bf Theorem}. {\stampatello (Null-function's Absolute Ramanujan Cloud is infinite)}
\par
\noindent
{\it The null-function} $\0:\N \rightarrow \C$ {\it has an infinite Absolute Ramanujan Cloud, i.e., } $\left|<\0>_{\rm abs}\right|=+\infty$. 
\smallskip
\par
\noindent{\bf Proof}. We may choose in previous Lemma $G=G_{\sigma}$ multiplicative, $G_{\sigma}(p_0^K)=1$, $\forall K\in \N_0$, for a {\stampatello fixed} prime $p_0$ and we have the Lemma's absolute convergence in the indices $(q,p_0)=1$, using $|c_q(a)|\le a$ $\forall a,q\in \N$, from $(3)$, choosing esp.  
$$
\forall \sigma>1\enspace \hbox{\stampatello fixed},
\quad
G_{\sigma}(p^K)=p^{-\sigma K},\enspace \forall K\in \N_0, \enspace \forall p\neq p_0\enspace \hbox{\stampatello prime}. 
$$ 
\par
\noindent
Thus, as $\sigma>1$ varies in the real numbers, we get an infinity of \enspace $G_{\sigma}\in <\0>_{\rm abs}$ \enspace from our Lemma.\hfill \square

\medskip

\par
We need to recall, here, two very well-known definitions, in number theory: the {\it kernel} of a natural number $n$ (also called the {\it radical} of $n$, compare [CG2] esp.), that is
$$
\kappa(n)\defineq \prod_{p|n}p,
\quad
\forall n\in \N,
$$
\par
\noindent
which has $\kappa(1)\defineq 1$, recalling that {\it empty products are $1$ by definition}; and the {\it square-full} naturals $s$ are defined by the property: 
$$
s\enspace {\it is}\enspace {\it square-full}
\quad \definiz \quad
v_p(s)\ge 2, \forall p|s \,; 
$$
\par
\noindent
in other words, for all $p\in \Primes$ fixed, $v_p(s)\neq 1$ (because it can only be $0$ or $\ge 2$): the number $1$ is square-full (all its $v_p$ vanish), strangely enough being also the only one square-free, too (all its $v_p$ are less than two)! For an arithmetic function $G$, in particular, the property of being \lq \lq {\it square-full supported}\rq \rq, then, simply means that: $q\in \N$ not square-full entails $G(q)=0$. 

\medskip

\par
We need the well-known result of Hildebrand [Hi] as proved in Th.m 1.1, Chapter V [ScSp], giving for any $F:\N \rightarrow \C$ an absolutely convergent (since it's finite) Ramanujan Expansion for $F$; that we state as: 
\smallskip
\par
\noindent{\bf Proposition}. {\stampatello (Hildebrand's Expansion for any Arithmetic Function)} 
\par
\noindent
{\stampatello Every arithmetic function} $F$ {\stampatello has a unique square-full supported Ramanujan Coefficient}, {\it say, } $\Hil \; F:\N \rightarrow \C$, {\stampatello its Hildebrand Coefficient}, {\it such that} : 
$$
\forall a\in \N
\qquad
F(a)=\sum_{q|a}\left(\Hil_{q\kappa(q)}\,F\right)c_{q\kappa(q)}(a). 
$$
\smallskip
\par
\noindent
In the Appendix, we give a Proof with complete details and also an Explicit Formula. 
\medskip
\par
We need this result, to prove the most important consequence of our Theorem.
\smallskip
\par
\noindent{\bf Corollary}. {\stampatello (Any arithmetic function has infinite Absolute Ramanujan Cloud)}
$$
F:\N \rightarrow \C
\quad \Longrightarrow \quad 
\left|<F>_{\rm abs}\right|=+\infty. 
$$

\vfill
\eject

\par				
\noindent{\bf Proof}. First of all, previous Proposition supplies, $\forall a\in \N$, {\stampatello an absolutely convergent Ramanujan Expansion for} $F(a)$, entailing 
$$
\Hil\; F\in <F>_{\rm abs} 
\quad \Rightarrow \quad
<F>_{\rm abs} \neq \emptyset. 
$$
\par
\noindent
Second, our Theorem ensures that we may choose in infinitely many ways $G\in <\0>_{\rm abs}$ ; notice that \enspace $\Hil\; F\in <F>_{\rm abs}$ {\stampatello and $G\in <\0>_{\rm abs}$ entail} $\forall a\in \N$, fixed, 
$$
\sum_{q=1}^{\infty}\left(\left(\Hil_q\,F\right)+G(q)\right)c_q(a)=F(a)
\enspace \hbox{\stampatello and} \enspace
\sum_{q=1}^{\infty}\left|\left(\Hil_q\,F\right)+G(q)\right|\,|c_q(a)|\le \sum_{q=1}^{\infty}\left(\left|\Hil_q\,F\right|+|G(q)|\right)|c_q(a)|<\infty, 
$$
\par
\noindent
whence : \enspace $\left(\left(\Hil_q\,F\right)+G(q)\right)\in <F>_{\rm abs}$ \enspace are an {\stampatello infinity of Ramanujan Coefficients for} $F$.\hfill $\square$

\vfill
\eject

\par				
\noindent{\bf 3. Remarks, enhancements, generalizations, glimpses and Ends} 
\bigskip
\par
\noindent
We start this complements corner, so to speak, with some Remarks.
\medskip
\par
\noindent{\bf Remark 1.} In the Lemma, the hypothesis \lq \lq $G$ multiplicative\rq \rq, alone, doesn't imply $G(1)=1$ (recall our definition); but, actually, the other hypothesis on $G$ does.\hfill $\diamond$
\medskip
\par
\noindent{\bf Remark 2.} We have an enhancement, for our Theorem, that extends naturally to our Corollary; namely, the cardinality of both $<\0>_{\rm abs}$ and, then, of our $<F>_{\rm abs}$ for general $F$ is uncountable: the coefficients $G_{\sigma}$ that we built (see Theorem Proof), actually, as our parameter $\sigma>1$ varies in real numbers, form an uncountable set. For any $F:\N\rightarrow \C$, of its Ramanujan Coefficients, see [C1] \& [CG2], there are plenty!\hfill $\diamond$
\medskip
\par
\noindent{\bf Remark 3.} Actually, we may also vary not only $\sigma$ in quoted Proof, but also the fixed $p_0$; however, prime numbers are countable. Notwithstanding this, we may consider more than one single $p_0$, say two, or more generally a finite set of such primes for $G_{\sigma}$, compare [CG1]; also, an infinite set of them and, even more, as we know the set of parts of countable sets is uncountable...\hfill $\diamond$
\medskip
\par
\noindent{\bf Remark 4.} We may also consider more general coefficients $G\in <\0>_{\rm abs}$, compare [CG1], giving rise to...\hfill $\diamond$
\smallskip
First Generalization, i.e.:
\smallskip
\par
\noindent{\bf Lemma}. {\stampatello (General Version)}
\par
\noindent
{\it Let } $G:\N \rightarrow \C$ {\it have the property that } $\exists p_0\in \Primes$, {\it with}
$$
G(mp_0^K)=G(m),
\qquad
\forall K\in \N_0,
\quad
\forall m\in \N, (m,p_0)=1.
$$
\par
\noindent
{\it Furthermore, let the series with indices coprime to $p_0$ , say}
$$
\sum_{{q=1}\atop {(q,p_0)=1}}^{\infty}G(q)c_q(a),
\qquad
\hbox{\it converge\enspace absolutely}
\enspace
\forall a\in \N\enspace \hbox{\it fixed}.
$$
\par
\noindent
{\it Then}
$$
G\in <\0>_{\rm abs}. 
$$
\smallskip
\par
\noindent
The {\stampatello Proof} is exactly the same as that of Lemma above ! From this, general Theorem \& Corollary follow. 
\medskip
\par
The first hypothesis on $G$, in this general version, is actually the definition of {\it weakly exotic}; while, first one with multiplicativity, of previous Lemma version, is that of \lq \lq {\it exotic}\rq \rq: for both, compare [CG1].
\medskip
\par
\noindent{\bf Remark 5.} From previous Remarks, it's clear that $<F>_{\rm abs}$, for any $F:\N\rightarrow \C$, is uncountable. However, from previous generalization, it's clearly linked to the cardinality of {\it weakly exotic coefficients} with the $p_0-$coprimality absolute convergence: are they only uncountable or even more?\hfill $\diamond$
\medskip
\par
\noindent{\bf Remark 6.} Previous Remark's question is for our future papers, on this interesting set (for general $F$) : the {\it Absolute Ramanujan Cloud} of our $F$. And, actually, Remark 5 is also our first glimpse on that. An initial attempt to describe and classify exactly the {\it Multiplicative Ramanujan Cloud} of our $F$, for general $F$, is the content of our previous paper [CG2].\hfill $\diamond$
\medskip
\par
\noindent{\bf Remark 7.} Back in time, we recall that Ramanujan's Coefficient for the null-function at paper's first page was also the first non-trivial one; some years later, Hardy [H] gave another. A question (as I don't know!), to our interested colleagues: any other more coefficients for $\0$, in the literature around those years?\hfill $\diamond$

\bigskip
\vfill

\par
In passing, we notice that Hildebrand's coefficient is a bit aside of general theories, for Ramanujan Expansions; one of the reasons is that given by Schwarz \& Spilker in their Book [ScSp] : it is not an expected one, as it doesn't fit, so to speak, neither Carmichael's [Ca] approach nor Wintner's [W] one; these two approaches, in fact, lead to, say, {\it Carmichael's coefficients} \& {\it Wintner's coefficients}, as we define them, esp., in [C1], where it is also clear that (apart from, say, \lq \lq pathologic arithmetic functions\rq \rq) they coincide, so leading to the so-called Carmichael-Wintner coefficients. (Compare our Appendix.) These are the coefficients that Schwarz \& Spilker take for granted as natural candidates to being Ramanujan Coefficients.

\eject

\par				
\noindent{\bf Appendix.} 
\par
\noindent
\bigskip
\par
\noindent$\underline{\hbox{\stampatello Proof of the Proposition}.}$ We prove: {\stampatello once fixed} $F:\N\rightarrow \C$,
\par
$\underline{\hbox{\stampatello if}}$ \enspace $G:\N\rightarrow \C$ \enspace {\stampatello is square-full supported and } $\forall a\in \N$, 
\enspace ${\displaystyle \sum_{q=1}^{\infty}G(q)c_q(a) }$ \enspace {\stampatello converges}, {\it w.r.t.} {\stampatello an arbitrarily fixed summation method, to } $F(a)$, 
\par
$\underline{\hbox{\stampatello then}}$ \qquad $\forall a\in \N$, \enspace ${\displaystyle F(a)=\sum_{m|a}G(m\kappa(m))c_{m\kappa(m)}(a) }$.
\par
\noindent
Of course, {\stampatello a posteriori the series is a finite sum}, converging absolutely, whence w.r.t. any summation method. So, we need another starting-point.
\par
{\stampatello The key-point is} the following
\smallskip
\par
\noindent$\underline{\hbox{\stampatello Fact}.}$ {\it Given a } {\stampatello square-full supported } $G:\N\rightarrow \C$, {\it then } $\forall a,q\in \N$ {\stampatello fixed}, 
$$
G(q)c_q(a)\neq 0
\quad \Rightarrow \quad
q=m\kappa(m) \enspace \hbox{\stampatello and} \enspace m|a.
\leqno{(\ast)}
$$
\par
\noindent
From this (in) fact, the above statement follows, whence the thesis. 
\par
{\stampatello We assume } $G(q)c_q(a)\neq 0$ {\stampatello henceforth}, when proving $(\ast)$ above. 
\par
\noindent
{\stampatello Actually}, from $G$ hypothesis and {\script R}amanujan {\script V}ertical {\script L}imit $(2)$, resp., the inequalities
$$
2\le v_p(q)\le v_p(a)+1,
\quad \forall p|q
$$
\par
\noindent
follow, with, resp., LHS $\le$ from hp. and RHS(Right Hand Side)$\le$ from {\script R}{\script V}{\script L}. Since $\kappa(q)|q$, $\forall q\in \N$, we write this as: 
$$
1\le v_p\left({q\over {\kappa(q)}}\right)\le v_p(a),
\quad \forall p|q;
$$
\par
\noindent
{\stampatello setting} ${q\over {\kappa(q)}}=m$, the equivalence \enspace $p|q\, \Leftrightarrow\,p|m$\enspace (from $q$ {\stampatello square-full}) {\stampatello entails}
$$
1\le v_p(m)\le v_p(a),\enspace \forall p|m
\quad \Rightarrow \quad m|a
$$
\par
\noindent
{\stampatello and by definition } $q=m\kappa(m)$, because \enspace $\kappa(m)=\kappa({q\over {\kappa(q)}})=\kappa(q)$, again from: $q$ {\stampatello square-full}.\hfill $\square$ 

\bigskip

\par
Once proved our {\stampatello Proposition}, we wish to define the {\stampatello Unique Hildebrand Coefficient}. (We follow Schwarz \& Spilker, quoted [ScSp].)
\par
This might seem even more difficult than previous Proof; but, like many times when we are in trouble in number theory, we dare to hope in the \lq \lq Mother of naturals\rq \rq: the {\stampatello Recursion}. 
\par
In fact, the only thing we know about Hildebrand $n-$th coefficient $\Hil_n\,F$, as square-full $n\in \N$ varies, once {\stampatello fixed } $F:\N\rightarrow \C$, is the formula 
$$
\forall a\in \N, 
\quad
F(a)=\sum_{q|a}\left(\Hil_{q\kappa(q)}\,F\right)c_{q\kappa(q)}(a). 
$$
\par
\noindent
Well, evidently,\enspace $F(1)=\left(\Hil_1\,F\right)c_1(1)$ \enspace $\Rightarrow $ \enspace $\Hil_1\,F\defineq F(1)$ \enspace is the \lq \lq {\stampatello Recursion starting point}\rq \rq, while 
$$
a>1
\enspace \Rightarrow \enspace
F(a)-\left(\Hil_{a\kappa(a)}\,F\right)c_{a\kappa(a)}(a)=\sum_{{q|a}\atop {q<a}}\left(\Hil_{q\kappa(q)}\,F\right)c_{q\kappa(q)}(a)
$$
\par
\noindent
entails the \lq \lq {\stampatello Recursion induction step}\rq \rq, say:
$$
a>1
\enspace \Rightarrow \enspace
\Hil_{a\kappa(a)}\,F \defineq {1\over {c_{a\kappa(a)}(a)}}\left(F(a)-\sum_{{q|a}\atop {q<a}}\left(\Hil_{q\kappa(q)}\,F\right)c_{q\kappa(q)}(a)\right).
$$
\par
\noindent
As [ScSp] underlines, we've used: $a>1$ $\Rightarrow$ $c_{a\kappa(a)}(a)\neq 0$, because, in fact, from $(1)$, we get
$$
c_{a\kappa(a)}(a)=\varphi(a\kappa(a)){{\mu(\kappa(a))}\over {\varphi(\kappa(a))}}\neq 0,
$$
\par
\noindent
by $\mu$ definition. 
\par
We explicitly point out that {\stampatello Hildebrand's Expansion for $F$}, in Proposition above, {\stampatello is a finite Ramanujan Expansion;} notwithstanding this, {\stampatello its length depends on} $a\in\N$ : all the properties enjoyed by fixed-length {\stampatello F.R.E.}s don't hold, for it; in fact, compare [C1], all {\stampatello fixed-length Ramanujan Expansions} of our $F$ {\stampatello have for Ramanujan Coefficients} the {\stampatello Carmichael-Wintner Coefficients for } $F$.  
\vfill
\eject

\par				
\centerline{\stampatello Bibliography}

\bigskip

\item{[Ca]} R. D. Carmichael, {\sl Expansions of arithmetical functions in infinite series}, Proc. London Math. Society {\bf 34} (1932), 1--26. 
\smallskip 
\item{[C1]} G. Coppola, {\sl Recent results on Ramanujan expansions with applications to correlations}, Rend. Sem. Mat. Univ. Pol. Torino {\bf 78.1} (2020), 57--82. 
\smallskip
\item{[C2]} G. Coppola, {\sl General elementary methods meeting elementary properties \thinspace of \thinspace correlations}, {\tt available}\break{\tt online at} \enspace arXiv:2309.17101 (2nd version)
\smallskip
\item{[C3]} G. Coppola, {\sl On Ramanujan smooth expansions for a general arithmetic function}, {\tt available online at} \enspace arXiv:2407.19759 
\smallskip
\item{[CG1]} G. Coppola and L. Ghidelli, {\sl Multiplicative Ramanujan coefficients of null-function}, arXiv:2005.14666v2 (2nd Version) 
\smallskip
\item{[CG2]} G. Coppola and L. Ghidelli, {\sl Convergence of Ramanujan expansions, I [Multiplicativity on Ramanujan clouds]}, arXiv:2009.14121v1 
\smallskip
\item{[CM]} G. Coppola and M. Ram Murty, {\sl Finite Ramanujan expansions and shifted convolution sums of arithmetical functions, II}, J. Number Theory {\bf 185} (2018), 16--47. 
\smallskip
\item{[D]} H. Davenport, {\sl Multiplicative Number Theory}, 3rd ed., GTM 74, Springer, New York, 2000. 
\smallskip
\item{[H]} G.H. Hardy, {\sl Note on Ramanujan's trigonometrical function $c_q(n)$ and certain series of arithmetical functions}, Proc. Cambridge Phil. Soc. {\bf 20} (1921), 263--271.
\smallskip
\item{[Hi]} A. Hildebrand, {\sl \"{U}ber die punktweise Konvergenz von Ramanujan-Entwicklungen zahlentheoretischer Funktionen}, Acta Arith. {\bf 44.2} (1984), 109--140.
\smallskip
\item{[L]} M. Laporta, {\sl On Ramanujan expansions with multiplicative coefficients}, Indian J. Pure Appl. Math. {\bf 52} (2021), 486--504. 
\smallskip
\item{[M]} M. Ram Murty, {\sl Ramanujan series for arithmetical functions}, Hardy-Ramanujan J. {\bf 36} (2013), 21--33. Available online 
\smallskip
\item{[R]} S. Ramanujan, {\sl On certain trigonometrical sums and their application to the theory of numbers}, Transactions Cambr. Phil. Soc. {\bf 22} (1918), 259--276.
\smallskip
\item{[ScSp]} W. Schwarz and J. Spilker, {\sl Arithmetical Functions}, Cambridge University Press, 1994.
\smallskip
\item{[T]} G. Tenenbaum, {\sl Introduction to Analytic and Probabilistic Number Theory}, Cambridge Studies in Advanced Mathematics, {46}, Cambridge University Press, 1995. 
\smallskip
\item{[W]} A. Wintner, {\sl Eratosthenian averages}, Waverly Press, Baltimore, MD, 1943. 

\bigskip
\bigskip
\bigskip

\par
\leftline{\tt Giovanni Coppola - Universit\`{a} degli Studi di Salerno (affiliation)}
\leftline{\tt Home address : Via Partenio 12 - 83100, Avellino (AV) - ITALY}
\leftline{\tt e-mail : giocop70@gmail.com}
\leftline{\tt e-page : www.giovannicoppola.name}
\leftline{\tt e-site : www.researchgate.net}

\bye